\documentclass[11pt]{article}
\usepackage{amsmath}
\usepackage{amssymb}
\usepackage{graphicx}
\usepackage{diagbox}
\usepackage{hyperref}
\usepackage{multirow}

\renewcommand{\baselinestretch}{1.1}
  \renewcommand{\arraystretch}{1.0}
\begin{document}

 \title{A Note on the Alon-Kleitman Argument\\ for Sum-free Subset Theorem}

  \author{Zhengjun Cao$^{1}$, \quad Lihua Liu$^{2,*}$}
  \footnotetext{ $^1$Department of Mathematics, Shanghai University, Shanghai, 200444,  China.  \\
   $^2$Department of Mathematics, Shanghai Maritime University, Shanghai, 201306,  China. \  $^*$\,\textsf{liulh@shmtu.edu.cn} }

 \date{}\maketitle

\begin{quotation}
 \textbf{Abstract}.
 In 1990, Alon and Kleitman proposed an argument for the sum-free subset problem:  every set of $n$ nonzero elements of a finite Abelian group contains a sum-free subset $A$ of size $|A|>\frac{2}{7}n$.  In this note, we  show that the argument confused two different randomness.  It applies only to the finite Abelian group  $G = (\mathbb{Z}/p\mathbb{Z})^s$ where $p$ is a prime.  For the general case,  the problem remains open.


  \end{quotation}

\section{Introduction}
 The sum-free subset problem is a special one in combinatorial number theory.   A subset $A$ of an Abelian group $G$ is called sum-free if
there are no   $a_1, a_2, a_3  \in  A$  such that $a_1 + a_2 = a_3$ ($a_1=a_2$ is not excluded).
In 1965, P. Erd\H{o}s \cite{E65} argued that every set of $n$ nonzero real numbers contains a sum-free subset $A$ of size $|A|\geq \frac{1}{3}n$.
 In 1990, Alon and Kleitman \cite{AK90} proposed their arguments for that:
 \begin{itemize}
  \item[(I)]{} every set of $n$ nonzero integers contains a sum-free subset $A$ of size $|A|>\frac{1}{3}n$;
   \item[(II)]{} every set of $n$ nonzero elements of a finite Abelian group contains a sum-free subset $A$ of size $|A|>\frac{2}{7}n$.
   \end{itemize}

Let $f(n)$ be the largest $\ell$ such that every set of $n$ nonzero elements contains
a sum-free subset of size $\ell$.  Alon and Kleitman \cite{AK90} pointed out that Erd\H{o}s's argument
can be modified to show $f(n) \geq \frac{1}{3}(n + 1)$.
For every set of $n$ positive integers, Bourgain \cite{B97} obtained  that  $f(n) \geq \frac{1}{3}(n + 2)$
 using an elaborate Fourier-analytic technique.  In 2014, Eberhard, Green and Manners \cite{EGM14} proved that
 there is a set of $n$ positive integers with no sum-free subset
of size greater than $\frac{1}{3}n + o(n)$ by constructing iteratively some approximate algebraic structures.

In 1971, Rhemtulla and Street \cite{RS71} obtained that if $G = (\mathbb{Z}/p\mathbb{Z})^s$, where $p=3k+1$ is a prime, then the
maximum sum free subset of $G$ has size $kp^{s-1}$. Taking $k = 2$, it gives $f(n)\leq \frac{2}{7}(n+1)$. In 1990,
 Alon and Kleitman \cite{AK90} claimed that every set $B$ of $n$ nonzero elements of a finite Abelian group contains a sum-free subset $A$ of size $|A|>\frac{2}{7}n$. The constant $\frac{2}{7}$ is best possible. The results on sum-free subset problem are summarized as follows (see Table 1).

\begin{center}
Table 1. Some results on sum-free subset problem \vspace*{3mm}

\begin{tabular}{|l|l|}
  \hline
     $n$ nonzero numbers &  Erd\H{o}s [1965], $f(n)\geq \frac{n}{3}$ \\ \cline{2-2}
     & Eberhard-Green-Manners [2014], $f(n) \leq (\frac{1}{3}+o(1)) n$ \\ \hline

  $n$ nonzero integers & Alon-Kleitman [1990],  $f(n) > \frac{1}{3}n $ \\   \hline
   $n$ positive integers & Bourgain [1997], $f(n) > \frac{1}{3}(n + 1)$  \\ \hline
   \hline
  $n$ nonzero elements in &  Rhemtulla-Street [1971], $f(n) \leq  \frac{2}{7}(n+1)$\\ \cline{2-2}
  a finite Abelian group &  Alon-Kleitman [1990],  $f(n) > \frac{2}{7}n$ \\
    \hline
\end{tabular}\end{center}

 Among these arguments, the Alon-Kleitman method is of special interest because it is elementary and can be used to design a deterministic algorithm \cite{AS08} for
finding such sum-free subsets. In this note, we would like to point out that Alon-Kleitman argument confuses two different randomness. The argument for the claim (II) 
applies only to a special case. For the general case, the problem remains open.

\section{Alon-Kleitman argument for a finite set of integers}
\subsection{Review}

\textbf{Theorem 2} [Alon and Kleitman, 1990] \textit{Every set $B=\{b_1, \cdots, b_n\}$ of $n$ nonzero integers contains a sum-free subset $A$ of size $|A|>\frac{1}{3}n$.}

\textit{Proof}. Let $p=3k+2$ be a prime, which satisfies $p>2\mbox{max}\{|b_i|\}_{i=1}^n$ and put
$C=\{k+1, k+2, \cdots, 2k+1\}.$ 
Observe that $C$ is a sum-free subset of the cyclic group $\mathbb{Z}_p$ and that
$\frac{|C|}{p-1}=\frac{k+1}{3k+1}>\frac{1}{3}. $ Choose at random an integer $x$, $1\leq x<p$, according to a uniform distribution on $\{1, 2, \cdots, p-1\}$, and define $d_1, \cdots, d_n$ by
$$d_i\equiv xb_i \ (\mbox{mod}\, p), 0\leq d_i<p. \eqno(1) $$
For every fixed $i$, $1\leq i\leq n$, as $x$ ranges over all numbers $1, 2, \cdots, p-1$, $d_i$ ranges over all nonzero elements of $\mathbb{Z}_p$ and hence
$$\mbox{Pr}[d_i\in C]=|C|/(p-1)>\frac{1}{3}. \eqno(2) $$
Therefore the expected number of elements $b_i$ such that $d_i\in C$ is more than $n/3$. Consequently, there is an $x$, $1\leq x<p$ and a subsequence $A$ of $B$ of cardinality $|A|>n/3,$ such that $xa \ (\mbox{mod}\, p)\in C$ for all $a\in A$. This $A$ is clearly sum-free, since if $a_1+a_2=a_3$ for some $a_1, a_2, a_3\in A$ then $xa_1+xa_2\equiv xa_3 \  (\mbox{mod}\, p)$, contradicting the fact that
$C$ is a sum-free subset of $\mathbb{Z}_p$. This completes the proof. \hfill $\Box$

\subsection{Analysis}

 It is easy to see that in
the definition  $d_i\equiv xb_i \ (\mbox{mod}\, p), 0\leq d_i<p, $
$b_i$ is viewed as the variable and $x$ is fixed. To specify the relationship, it is better to
write it as
$$d_x(b_i)\equiv x\,b_i \ (\mbox{mod}\, p)   \eqno(3) $$
If $x$ is viewed as the variable and $b_i$ is fixed, then it should define
$$d_{b_i}(x)\equiv x\,b_i \ (\mbox{mod}\, p)  \eqno(4) $$

Trivially, as $x$ ranges over all numbers $1, 2, \cdots, p-1$, $d_i(x)$ ranges over all nonzero elements of $\mathbb{Z}_p$ (see the definitions of $p$).
Define
$B_i:=\{b_i\ (\mbox{mod}\, p), 2b_i\ (\mbox{mod}\, p), \cdots, (p-1)b_i\ (\mbox{mod}\, p) \}.$
We have $B_i=\{1, 2, \cdots, p-1\}, $
$$\mbox{Pr}[d_{b_i}(x)\in C]=|B_i\cap C|/|B_i|=|C|/(p-1)>\frac{1}{3}. \eqno(5) $$
Define
$A_x:=\{ xb_1\ (\mbox{mod}\, p), xb_2\ (\mbox{mod}\, p), \cdots, xb_n\ (\mbox{mod}\, p) \}$ for some picked $x\in \mathbb{Z}_p^{*}$.
Clearly, $|A_x|=n$ and
$$\mbox{Pr}[d_x(b_i)\in C]=|A_x\cap C|/|A_x| =|A_x\cap C|/n. \eqno(6) $$
One cannot conclude that
$\mbox{Pr}[d_x(b_i)\in C]>\frac{1}{3},$
because $|A_x\cap C|$ is still undetermined.
Frankly speaking,  the argument confuses the trivial probability  $\mbox{Pr}[d_{b_i}(x)\in C]$ with the wanted probability $\mbox{Pr}[d_x(b_i)\in C]$.  

\subsection{An explicit proof}

\textit{Proof}. Let $p=3k+2$ be a prime, which satisfies $p>2\mbox{max}\{|b_i|\}_{i=1}^n$ and put
$C=\{k+1, k+2, \cdots, 2k+1\}.$
Construct the following residues table (see Table 2). It is easy to find that each row has exactly $k+1$ residues in $C$. Totally, the table has $n(k+1)$ residues in $C$.

\begin{center}
Table 2. Residues table for integers\vspace*{3mm}

\begin{tabular}{|l|cccc|c|}
\hline
     \diagbox{$b$}{residue}{$x$}    & 1   & 2 & $\cdots$ & p-1&  \\ \hline
   $b_1$ & $b_1 (\mbox{mod}\, p)$ & $2b_1(\mbox{mod}\, p)$   &          & $(p-1)b_1 (\mbox{mod}\, p)$ &  $k+1$ residues in $C$\\
   $b_2$ & $b_2 (\mbox{mod}\, p)$ & $2b_2(\mbox{mod}\, p)$  &          & $(p-1)b_2 (\mbox{mod}\, p)$ &  $k+1$ residues in $C$\\
   $\vdots$ &   &   &         &  & \\
   $b_n$ & $b_n (\mbox{mod}\, p)$ & $2b_n(\mbox{mod}\, p)$  &          & $(p-1)b_n (\mbox{mod}\, p)$ &   $k+1$ residues in $C$\\ \hline
     &    &    &   $\cdots$      &   & $n(k+1)$ residues in $C$\\
    \hline
\end{tabular}\end{center}

Suppose that each column has at most $\frac n 3$ residues in $C$.
Since $b_i x \not\equiv b_j x (\mbox{mod}\, p)$ for $  x  \in \mathbb{Z}_p^*$, if $i, j \in \{1, 2, \cdots, n\}$ and $i \neq j$,
it follows that the table has at most $\frac n 3 (p-1)$ residues in $C$. Hence,
$n(k+1)\leq \frac n 3 (p-1) $, $3k+4\leq p$. It contradicts that $p=3k+2$. Thus,  there is a column
which has strictly at least $\frac n 3$ residues in $C$. Denote the column by
$$ X:=\{\hat{x} b_1 (\mbox{mod}\, p), \hat{x} b_2 (\mbox{mod}\, p),\cdots, \hat{x} b_n (\mbox{mod}\, p) \}.  $$
We have $|X\cap C|>\frac n 3$. Denote the intersection of $X$ and $C$ by
$$\hat X:=\{\hat{x} b_{i1} (\mbox{mod}\, p), \hat{x} b_{i2} (\mbox{mod}\, p),\cdots, \hat{x} b_{ik} (\mbox{mod}\, p) \}. $$
We obtain the set $A:=\{b_{i1}, b_{i2}, \cdots, b_{ik} \}$.
Since $\hat X \subset C$, $C$  is sum-free and the mapping
\begin{eqnarray*}
\mathcal{M}_p: \quad   A &\rightarrow& \hat X\\
  b_{ij} &\mapsto&  \hat{x} b_{ij} (\mbox{mod}\, p)
  \end{eqnarray*}
  is \underline{a bijection}, we conclude that  $A$ is sum-free. Clearly, $|A|=|\hat X|>\frac n 3.$ \hfill $\Box$

\section{Alon-Kleitman argument for a finite Abelian group}
\subsection{Review}

For a subset $B$ of an Abelian group, let $s(B)$ denote the maximum cardinality of a sum-free subset of  $B$. Similarly, for a sequence $A=(a_1, a_2, \cdots, a_n)$ of (not necessarily distinct) elements of an Abelian group, let $s(A)$ denote the maximum number of elements in a sum-free subsequence $(a_{i_1}, a_{i_2}, \cdots, a_{i_k} )$ of $A$.
 Define
 \begin{eqnarray*}
 I_1&=&\{x\in \mathbb{Z}_n: \frac 1 3 n<x\leq \frac 2 3 n\},\\
 I_2&=&\{x\in \mathbb{Z}_n: \frac 1 6 n<x\leq \frac 1 3 n\ \mbox{or}\ \frac 2 3 n<x\leq \frac 5 6 n \}.
 \end{eqnarray*}
 It is easy to check that both $I_1$ and $I_2$ are sum-free subsets of $\mathbb{Z}_n$.
For any divisor $d$ of $n$, let $d\,\mathbb{Z}_n$ denote the subgroup of all multiples of $d$ in $\mathbb{Z}_n$, i.e., $d\,\mathbb{Z}_n=\{0, d, 2d, \cdots, n-d\}$. Clearly $d\,\mathbb{Z}_n$ has $n/d$ elements. By a straightforward case analysis, for all admissible values of $n$ and $d$, we have
$$\frac 4 7 \frac{|d\mathbb{Z}_n\cap I_1|}{|d\mathbb{Z}_n|}+ \frac 3 7 \frac{|d\mathbb{Z}_n\cap I_2|}{|d\mathbb{Z}_n|} \geq \frac 2 7.$$
%
%

\textbf{Theorem 3} [Alon and Kleitman, 1990] \textit{For any finite Abelian group $G$, every set $B$ of non-zero elements of $G$ satisfies $s(B)>\frac 2 7 |B|$. The constant $\frac 2 7$ is best  possible. Similarly, every sequence $A$ of non-zero elements of $G$ satisfies $s(A)>\frac 2 7 |A|$, and the constant $\frac 2 7$ is optimal.}

 \emph{Proof}.  It suffices to show that, for any finite Abelian group $G$ and every sequence $A$ of non-zero elements of $G$, $s(A)>\frac 2 7 |A|$. 
Let $G$ be an arbitrary finite Abelian group and let $B=(\textbf{b}_1, \textbf{b}_2, \cdots, \textbf{b}_m)$ be a sequence of $m$ non-zero elements of $G$. As is well known, $G$ is a direct sum of cyclic groups and therefore there are $n$ and $s$ such that $G$ is a subgroup of the direct sum of $s$ copies of $\mathbb{Z}_n$. Thus we can think of the elements of $B$ as members of $\mathbb{Z}_n^s$. Each such element $\textbf{b}_i$ is a vector, i.e., $\textbf{b}_i=(b_{i1}, b_{i2}, \cdots, b_{is})$, where, for each $i$, $0\leq b_{i1}, \cdots, b_{is}<n$ and not all the $b_{ij}$ are zero.

Choose a random element $\textbf{x}=(x_1, x_2, \cdots, x_s)$ of $\mathbb{Z}_n^s$  and define  $f_1, f_2, \cdots, f_m$ by
$$f_i=\sum_{j=1}^s x_jb_{ij}\ (\mbox{mod}\, n).  \eqno(7)$$
Notice that for every fixed $i\, (1\leq i\leq m)$, the mapping
$$(x_1, x_2, \cdots, x_s)  \mapsto \sum_{j=1}^s x_jb_{ij} \ (\mbox{mod}\, n) \eqno(8)$$
is a homomorphism from $\mathbb{Z}_n^s$ to $\mathbb{Z}_n$. Moreover, if $d_i=\mbox{gcd}(b_{i1}, b_{i2}, \cdots, b_{is}, n)$ then the image of this homomorphism is just $d_i\mathbb{Z}_n$. Consequently, as $\textbf{x}=(x_1, \cdots, x_s)$   ranges over all
 elements of $\mathbb{Z}_n^s$, $f_i$ ranges over all elements of $d_i\mathbb{Z}_n$ and attains each value of $d_i\mathbb{Z}_n$  the same number of times. It follows that, for each $j=1, 2$,
 $$\mbox{Pr}(f_i\in I_j)=\frac{|d_i\mathbb{Z}_n\cap I_j|}{|d_i\mathbb{Z}_n|}. \eqno(9) $$

For each divisor $d$ of $n\,(1\leq d<n),$ let $m_d$ denote the number of elements $\textbf{b}_i$ in $B$ such that $\mbox{gcd}(b_{i1}, b_{i2}, \cdots, b_{is},n)=d.$ Clearly,
$\sum_{d|n} m_d=m.$   For $j=1, 2$,
 the expected number of elements $\textbf{b}_i$ such that $f_i\in I_j$ is
$$M_j=\sum_{d\,|\,n} m_d \mbox{Pr}(f_i\in I_j)=\sum_{d\,|\,n} m_d \frac{|d\,\mathbb{Z}_n\cap I_j|}{|d\,\mathbb{Z}_n|}. \eqno(10)$$
Moreover, since $(0,0,\cdots, 0)\in \mathbb{Z}_n^s$ maps every $b_i$ into $f_i=0\not\in I_1$, it follows that there is an $\textbf{x}=(x_1, x_2, \cdots, x_s)\in \mathbb{Z}_n^s$ and a subsequence $A$ of strictly more than $M_1$ elements of $B$ such that each $\textbf{a}=(a_1, \cdots, a_s)\in A$  is mapped by $\textbf{x}$ into $\sum_{i=1}^s x_ia_i (\mbox{mod}\, n) \in I_1.$ Clearly {this $A$ is a sum-free} subsequence of $B$ (since $I_1$ is sum-free) and thus
$s(B)>M_1=\sum_{d\,|\,n} m_d \frac{|d\,\mathbb{Z}_n\cap I_1|}{|d\,\mathbb{Z}_n|}.$
Similarly, since $I_2$ is sum-free,
$s(B)>M_2=\sum_{d\,|\,n} m_d \frac{|d\,\mathbb{Z}_n\cap I_2|}{|d\,\mathbb{Z}_n|}.$
Hence, we have
\begin{eqnarray*}
 s(B)=\frac 4 7 s(B) +\frac 3 7 s(B) &>& \frac 4 7 M_1+\frac 3 7 M_2\\
 &\geq & \sum_{d\,|\,n} m_d \left( \frac 4 7 \frac{|d\,\mathbb{Z}_n\cap I_1|}{|d\,\mathbb{Z}_n|} +\frac 3 7 \frac{|d\,\mathbb{Z}_n\cap I_2|}{|d\,\mathbb{Z}_n|} \right)\\
  &\geq & \sum_{d\,|\,n} m_d \times \frac 2 7=\frac 2 7 m.
 \end{eqnarray*}

\subsection{It confuses two different randomness} 

It is easy to see that in
the definition  $f_i=\sum_{j=1}^s x_jb_{ij}\ (\mbox{mod}\, n) $,
$\textbf{b}_i$ is viewed as the variable and $\textbf{x}$ is fixed. To specify the relationship, it is better to
write it as
$$f_{\textbf{x}}(\textbf{b}_i)=\sum_{j=1}^s x_jb_{ij}\ (\mbox{mod}\, n)   \eqno(7') $$
If $\textbf{x}$ is viewed as the variable and $\textbf{b}_i$ is fixed, then it should define
$$f_{\textbf{b}_i}(\textbf{x})=\sum_{j=1}^s x_jb_{ij}\ (\mbox{mod}\, n)  \eqno(8') $$
Consequently, the Eq.(9) should be corrected as
$$\mbox{Pr}(f_{\textbf{b}_i}(\textbf{x})\in I_j)=\frac{|d_i\mathbb{Z}_n\cap I_j|}{|d_i\mathbb{Z}_n|}  \eqno(9') $$
For $j=1, 2$,
 the expected number of elements $\textbf{b}_i$ such that $f_{\textbf{x}}(\textbf{b}_i)\in I_j$ is
$$M_j=\sum_{d\,|\,n} m_d \mbox{Pr}(f_{\textbf{x}}(\textbf{b}_i)\in I_j)\neq \sum_{d\,|\,n} m_d \frac{|d\,\mathbb{Z}_n\cap I_j|}{|d\,\mathbb{Z}_n|}. \eqno(10')$$
Apparently, in Alon-Kleitman argument it confuses the trivial probability  $\mbox{Pr}(f_{\textbf{b}_i}(\textbf{x})\in I_j)$ with the wanted probability $\mbox{Pr}(f_{\textbf{x}}(\textbf{b}_i)\in I_j)$. 

 \subsection{It is short of a bijection}

In the  argument for a finite Abelian group, for $\textbf{x}=(x_1, x_2, \cdots, x_s)\in  \mathbb{Z}_n^s$  the  mapping $\mathcal{M}_{\textbf{x}}$,
\begin{eqnarray*}
\mathcal{M}_\textbf{x}: \qquad \qquad\qquad  A &\rightarrow& \hat A, \ \mbox{where}\, \hat A\subset \mathbb{Z}_n \\
(a_1, \cdots, a_s)\in  \mathbb{Z}_n^s &\mapsto& \sum_{i=1}^s x_ia_i\, (\mbox{mod}\, n)
\end{eqnarray*}
is a homomorphism, \underline{not a bijection}.
 It is easy to find that the following Table 3 is far different from Table 2, where
$\textbf{x}_i=(x_{i1}, x_{i2}, \cdots, x_{is})\in \mathbb{Z}_n^s, i=1, \cdots, n^s-1,
\textbf{b}_j=(b_{j1}, b_{j2}, \cdots, b_{js})\in \mathbb{Z}_n^s, j=1, \cdots, m,$ $d_i=\mbox{gcd}\,(b_{i1}, b_{i2}, \cdots, b_{is}, n)$, $|I_1|=[\frac{2}{3}n]-[\frac{1}{3}n]$.

\begin{center}
Table 3. Residues table for vectors \vspace*{3mm}

\begin{tabular}{|c|cccc|l|}
\hline
     \diagbox{$\textbf{b}$}{residue}{$\textbf{x}$}  & $\textbf{x}_1$ & $\cdots$ & $\textbf{x}_i$ & $\cdots$ &  \\ \hline
   $\textbf{b}_1$           & $\textbf{b}_1\cdot \textbf{x}_1 \ (\mbox{mod}\, n)$      & $\cdots$  & $\textbf{b}_1\cdot \textbf{x}_i \ (\mbox{mod}\, n)$  & $\cdots$ & $d_1n^{s-1}|d_1\mathbb{Z}_n\cap I_1|$ residues in $I_1$ \\
    $\vdots$     &    &    &   & &  \\
   $\textbf{b}_m$            & $\textbf{b}_m \cdot \textbf{x}_1 \ (\mbox{mod}\, n)$   & $\cdots$  & $\textbf{b}_m \cdot \textbf{x}_i \ (\mbox{mod}\, n)$ & $\cdots$ & $d_mn^{s-1}|d_m\mathbb{Z}_n\cap I_1|$ residues in $I_1$\\ \hline
              &    &   &   &  & $\sum_{i=1}^{m}d_i n^{s-1}|d_i\mathbb{Z}_n\cap I_1|$ residues in $I_1$ \\
    \hline
\end{tabular}\end{center}

Since there are $n^s-1$ columns, for each column  the expected  number of residues in $I_1$  is
$ \frac{\sum_{i=1}^{m} d_in^{s-1}|d_i\mathbb{Z}_n\cap I_1|}{n^s-1}$.
Suppose that $\alpha=\mbox{min}\,\{d_1, d_2, \cdots, d_m\}, \beta=\mbox{max}\,\{d_1, d_2, \cdots, d_m\}$, we have
\begin{eqnarray*}
\frac{\sum_{i=1}^{m}d_i n^{s-1}|d_i\mathbb{Z}_n\cap I_1|}{n^s-1}&\geq& \frac{\alpha \sum_{i=1}^{m} n^{s-1}|d_i\mathbb{Z}_n\cap I_1|}{n^s-1}\\
&\geq& \frac{\alpha \sum_{i=1}^{m} n^{s-1}|\beta \mathbb{Z}_n\cap I_1|}{n^s-1}\\
&=& \alpha m \times\left[\frac{[\frac{2}{3}n]-[\frac{1}{3}n]}{\beta}\right]\times \frac{n^{s-1}}{n^s-1}\\
&\rightarrow& \frac{\alpha}{3\beta} m \ (n\rightarrow \infty).
\end{eqnarray*}
Hence, using the residues table one \underline{cannot obtain} the wanted equality
$$\frac{\sum_{i=1}^{m}d_i n^{s-1}|d_i\mathbb{Z}_n\cap I_1|}{n^s-1}\geq \frac{2}{7}m. $$
Thus, Alon-Kleitman argument method for a finite Abelian group fails.

In the above analysis, if $n$ is a prime then $d_i\mathbb{Z}_n=\mathbb{Z}_n$, $i=1, 2, \cdots, m$. It follows that
 for each column  the expected  number of residues in $I_1$  is
\begin{eqnarray*}
& & \frac{\sum_{i=1}^{m} d_in^{s-1}|d_i\mathbb{Z}_n\cap I_1|}{n^s-1}=\frac{\sum_{i=1}^{m} n^{s-1}|\mathbb{Z}_n\cap I_1|}{n^s-1}=\frac{m n^{s-1}| I_1|}{n^s-1}\\
& &=\frac{ 1}{1-n^{-s}}\times \frac{[\frac{2}{3}n]-[\frac{1}{3}n]}{n} \times m > \frac{[\frac{2}{3}n]-[\frac{1}{3}n]}{n} \times m.
\end{eqnarray*}
Since $n$ is a prime, then $\frac{[\frac{2}{3}n]-[\frac{1}{3}n]}{n}\geq \frac{2}{7}.$
Therefore, we can obtain the following result.

\textbf{Theorem 4}  \textit{For any finite Abelian group $G$, if $G = (\mathbb{Z}/p\mathbb{Z})^s$ where $p$ is a prime, then every set $B$ of non-zero elements of $G$ satisfies $s(B)>\frac 2 7 |B|$. The constant $\frac 2 7$ is best  possible. Similarly, every sequence $A$ of non-zero elements of $G$ satisfies $s(A)>\frac 2 7 |A|$, and the constant $\frac 2 7$ is optimal.}

\section{Conclusion}
W show that  Alon-Kleitman argument method for sum-free subset problem  is flawed.  We want to point out that, the claim that every set of $n$ nonzero elements of a finite Abelian group contains a sum-free subset $A$ of size $|A|>\frac{2}{7}n$, remains open.


\end{document}